\documentclass{amsart}
\usepackage{amsmath,amsthm,amsfonts,amssymb}
\def\({\left(}
\def\){\right)}

\def\Br{\mathbb R}
\def\Bz{\mathbb Z}

\def\flx{((X,T)_g,\Br)}
\def\fx{((X,T)_h,\Br)}
\def\fly{((Y,S)_h,\Br)}
\def\fy{((Y,S)_g,\Br)}

\theoremstyle{plain}   
\begingroup 
\newtheorem{thm}{Theorem}[section]   %
\newtheorem{prop}[thm]{Proposition}  
\endgroup

\theoremstyle{definition}

\theoremstyle{remark}
\newtheorem{rem}[thm]{Remark}        %

\numberwithin{equation}{section}

\begin{document}


\title[Factor maps between tiling dynamical systems]
{Factor maps between tiling dynamical systems\\
Draft---{\today}}

\author{Karl Petersen}

\address{Department of Mathematics,
CB 3250, Phillips Hall,
         University of North Carolina,
Chapel Hill, NC 27599 USA}
\email{petersen@math.unc.edu}





\begin{abstract}
We show that there is no Curtis-Hedlund-Lyndon Theorem for factor maps
between tiling dynamical systems: there are codes between such systems
which cannot be achieved by working within a finite window. By
considering 1-dimensional tiling systems, which are the same as flows
under functions on subshifts with finite alphabets of symbols, we
construct a `simple' code which is not `local', a local code which is
not simple, and a continuous code which is neither local nor simple.
\end{abstract}

\maketitle

\section{Introduction}\label{intro}
According to the Curtis-Hedlund-Lyndon Theorem \cite{CHL}, every
factor map (continuous shift-commuting map) between subshifts (closed
shift-invariant subsets of spaces of one- or two-sided sequences on
finite alphabets) is a sliding block code, that is to say the central
coordinate of the image of any point is determined by a finite range
(of fixed size) of central coordinates of the point. Tiling dynamical
systems also are based on a finite set of symbols, the prototiles. A
tiling is a covering of $\mathbb R^d$ by congruent copies, called
tiles, of members of the finite set of prototiles which intersect only
along their boundaries. Two tilings are considered to be close if
within a large neighborhood of the origin the unions of the tile
boundaries are close in Hausdorff distance. Then $\mathbb R^d$ acts 
continuously on
the space of tilings by translations, and tiling dynamical systems are
defined to be closed invariant subsets of the set of all tilings of
$\mathbb R^d$. Often continuous translation-commuting maps between
tiling dynamical systems are constructed as {\em local codes}: in the
image of a tiling, the tile type found at the origin, and the precise
location of the origin within that tile, are determined by a
fixed-radius neighborhood of the origin in the tiling in the
domain. See
\cite{Priebe,Radin,Robinson,Senechal,Solomyak} and the references
cited in those sources for background and for examples of tiling dynamical
systems and local codes.

Our purpose here is to show that there are (continuous,
translation-commuting) factor maps between tiling dynamical systems
which are not given by local codes. In fact, such examples can be
found for  1-dimensional tiling dynamical systems, which are flows under
functions built on subshifts. Mappings between flows under functions
commute with translations if and only if they satisfy some
cohomological equations (see (\ref{cohom}) below). Some particularly
simple examples of this kind (see (\ref{simple})) will be called {\em
simple} maps. We will make first an example of a simple map which is
not local, then a local map which is not simple, and finally a factor
map which is neither local nor simple. The examples we construct are
based on specific symbolic dynamical systems with particular ceiling
functions, but only certain properties of the
systems and functions are really necessary, so many similar examples
can be constructed easily.

I thank N.~Priebe and E.~A.~Robinson, Jr., for bringing this problem to
my attention, the Erwin Schr\" odinger Institute, Vienna, for
support while part of this research was conducted, and C.~Radin for
pointing out that the preprint \cite{Radin-Sadun} also gives an example of a 
non-local code.

\section{Tiling systems and factor maps}\label{systems}
Let $(X,T)$ and $(Y,S)$ be subshifts on finite alphabets $A$ and
$B$. Let $g:X \to (0,\infty )$ and $h:Y \to (0,\infty )$ be continuous
functions which depend only on the central entry: $g(x)=g_0(x_0),
h(y)=h_0(y_0)$.  (If we were given maps $g$ and $h$ that depended on
finitely many entries, we could form higher block representations of
$X$ and $Y$ to get the dependence down to just the central entry.)  We
denote by $((X,T)_g,\mathbb R)$ and $((Y,S)_h,\mathbb R)$ the flows
built under the ceiling functions $g$ and $h$. Recall that $(X,T)_g$,
for example, is the quotient space of $(X \times \mathbb R, \mathbb
R)$ (with the action $(x,s)t=(x,s+t)$) under the equivalence $\sim$ generated
by $(x,g(x)) \equiv (Tx,0)$.  It is natural to use the notation
$[x,s]$ for the equivalence class of a pair $(x,s)$ when $ x \in X, 0
\leq s < g(x)$, as in \cite{KMS1,KMS2}.  For each equivalence class
$\xi \in (X,T)_g$, there are a unique symbolic sequence $\pi _X\xi \in
X$ and a unique $\pi _{\Br}\xi \geq 0$ such that $0 \leq \pi
_{\Br}\xi < g(\pi _X\xi )$ and $\xi = [\pi _X\xi ,\pi
_{\Br}\xi ]$.
\begin{rem}
The existence
of these maps is special for the one-dimensional situation; in
higher-dimensional cases, except for very regular tilings 
we would obtain
a labeled graph rather than
a symbolic sequence, and perhaps barycentric coordinates of the origin
in the central tile rather than the coordinate $\pi _{\Br}\xi$.
 The maps $\pi
_X$ and $\pi_{\Br}$ are usually not continuous nor well-behaved with
respect to the actions.
More properly the notation for these maps would also display their domain,
but we rely on the context.
\end{rem}

With each point $[x,s] \in (X,T)_g$, we associate a tiling of $\mathbb
R$: we have, for each element $a$ of the alphabet $A$, a prototile
which is a closed interval of length $g_0(a)$.  The tiling
corresponding to $[x,s]$ $(x \in X, 0 \leq s < g(x))$ consists of the
sequence of tiles specified by the symbolic sequence $x$; and the
central tile, of type $x_0$, contains the origin $s$ units from its
left endpoint.  This identification is a topological conjugacy between
the flow under a function and the associated tiling dynamical system.

  A {\em factor map} or {\em code} $\phi : \flx \to \fly$ is a
continuous onto map which commutes with the action of $\Br$.  We will
say that a factor map $\phi$ is {\em local} if there is $r \geq 0$
such that if $x,x' \in X$ , $x_i=x'_i$ for $|i| \leq r$, $0 \leq s <
g(x)(=g(x')$), and $\phi [x,s]=[y,t]$ (for some $y \in Y$ and $0 \leq t
< h(y)$), then also $\phi[x',s]=[y,t]$.  Easy examples of local codes
arise from subdividing tiles or amalgamating patches of tiles into new
tiles. We say that the factor map
$\phi$ is {\em simple} (cf. \cite{KMS1,KMS2}) if there are a factor map
(continuous shift-commuting map) $\pi:(X,T) \to (Y,S)$ and a
function $t:X \to \Br$ such that
\begin{equation}\label{simple}
\begin{split}
t(Tx) - t(x) &= g(x) - h(\pi x) \quad \text{ for all } x \in X\quad\text{and}\\
\phi[x,s]&=[\pi x,0](s + t(x)) \quad \text{  for all } [x,s] \in \flx .
\end{split}
\end{equation}

  Given $g:X \to \Br$, define a cocycle $g(x,n)$ on $(X,T)$ by
\begin{equation}
 g(x,n)=\begin{cases} \sum_{k=0}^{n-1}g(T^kx) &\text{  if } n>0,\\
        0 &\text{  if } n=0 ,\\
        -\sum_{k=n}^{-1}g(T^kx) &\text{  if } n< 0 .
\end{cases}      
\end{equation}
If $g:X \to \Br$ is fixed, for $x \in X$ and $s \in \Br$ define
$n_X(x,s)=n_X^{(g)}(x,s)$ by
\begin{equation}
g(x,n_X(x,s))\leq s < g(x,n_X(x,s)+1) . 
\end{equation}
Then the action of $\Br$ on $\flx$ is given by
\begin{equation}
[x,s]u=\left[T^{n_X(x,s+u)}x,s+u-g(x,n_X(x,s+u))\right] \quad (x \in X,
0 \leq s < g(x))
\end{equation}
(note that $0 \leq s+u-g(x,n_X(x,s+u)) < g(T^{n_X(x,s+u)}x))$.
Similar definitions and formulas apply to $\fly$.
Further, in $\flx$ we have 
\begin{equation}
\begin{split}
[x,s]=[x',s'] \quad\text{  if and only if}\\
\text{ there is } n \in \Bz \text{
with
 } T^nx=x' \text{  and  } s'=s-g(x,n) .
\end{split}
\end{equation}

\begin{prop}
Given functions $\pi : X \to Y$ and $v: X \to \Br$ such
that 
\begin{equation}\label{cohom}
\begin{split}
S^{n_Y(\pi x,g(x)+v(x))}\pi(x) &= S^{n_Y(\pi Tx, v(Tx))}(\pi Tx)\quad\text{and}\\
v(Tx)-v(x)&=g(x)+h(\pi Tx,n_Y(\pi Tx,v(Tx)))\\
          &\quad\qquad\qquad -h(\pi x, n_Y(\pi x,g(x)+v(x))) ,
\end{split}
\end{equation}
then putting
\begin{equation}\label{factor}
\phi [x,s] = [\pi x, 0] (s+v(x))
\end{equation}
determines a map $\phi: \flx \to \fly$
which commutes with the $\Br$ actions and thus is a factor map if it
is onto and continuous (which will be the case if $\pi$ and $v$ are continuous).
Conversely, given a factor map $\phi: \flx \to \fly$, there are
functions  $\pi : X \to Y$ and $v: X \to \Br$ such
that (\ref{cohom}) holds for all $[x,s] \in \flx$.
\end{prop}

\begin{proof}
The equations (\ref{cohom}) guarantee that the map defined by
(\ref{factor}) respects equivalence
classes, so $\phi$ is well-defined from $\flx$ to $\fly$ and commutes
with the actions.
(Clearly $(\phi [x,s])u=\phi([x,s]u)$ as long as $0\leq s+u <
g(x)$. The equations (\ref{cohom}) imply that the commuting is
maintained when we move across the boundary of a tile in $\flx$.)
 
 Conversely, given a factor map $\phi: \flx \to \fly$, we
define
\begin{equation*}
\pi x = \pi_Y(\phi[x,0]), \quad v(x)=\pi _{\Br}(\phi [x,0]) .
\end{equation*}
Then
\begin{equation}
\begin{split}
\phi [x,s] &= (\phi [x,0])s = [\pi x,v(x)]s \quad \text{(with }0 \leq
v(x) <h(\pi x))\\
&=[\pi x,0](v(x)+s),
\end{split}
\end{equation}
and then $\pi$ and $v$ necessarily satisfy (\ref{cohom}).
\end{proof}

In the case of a simple map, $\pi: (X,T) \to (Y,S)$ is a factor map
and we take $v(x)=t(x)$.
Then Equations (\ref{cohom}) follow from (\ref{simple}).

\section{Examples}\label{examples}

\subsection{Example 1}\label{Ex1} 
We first construct an example of a simple map 
which is not a local code.

Let $(X,T)$ be a Sturmian subshift which codes translation
$R_{\alpha}(t)=\langle t+\alpha \rangle = t + \alpha \mod 1$ by an
irrational $\alpha$ on $[0,1]$: define
$\omega(n)=\chi_{[0,1/2)}\langle n\alpha\rangle$ for all $n \in \Bz$
and let $X$ be the orbit closure of $\omega$ under the shift
transformation $T$. Then $(X,T)$ is a minimal, uniquely ergodic
topological dynamical system, and there is a factor map $\rho :(X,T)
\to ([0,1),R_{\alpha})$ which is one-to-one except on a countable set
on which it is two-to-one (the union of the orbits of 0 and $1/2$
under $R_{\alpha}$).

 Let $(Y,S)=(X,T)$ and let the factor map $\pi :(X,T) \to (Y,S)$ be
the identity.  We will specify $h$ and $t$ first and then {\em define}
$g$ so that Equation (\ref{simple}) is satisfied.  Fix $\gamma \in \Bz
\alpha \cap (0,1/4)$ and $\eta_1, \eta_2 > 5$ and define $h$ by 
\begin{equation}
h(y)=\begin{cases}
	\eta_1\quad&\text{ if } \rho y \in {[0,\gamma )}\\
	\eta_2 &\text{ if } \rho y \in [\gamma,1) . 
	\end{cases} 
\end{equation} 
Note that $h$ is continuous on $Y$ and takes only finitely many
values, so by recoding we may assume that $h$ is a function of just
the central coordinate. Let 
\begin{equation} t(x) = \rho (x) , 
\end{equation} 
so that 
\begin{equation} t(Tx) -t(x) = \begin{cases} 
                  \alpha &\text{ if } 0 \leq \rho (x) < 1-\alpha\\
                   -1+\alpha &\text{ if } 1-\alpha \leq \rho (x) < 1 .
                   \end{cases} 
\end{equation}
Again $t$ is continuous on $X$, even though the function $f(s)=s$ is
not continuous on $S^1=$the reals modulo $1$.
(If $x^{(n)} \to x$ in $X$, then eventually $(x^{(n)})_0$ stabilizes at
either
$0$ or $1$, so that for large enough $n$ either all $\rho (x^{(n)})
\in [1/2,1)$ or $\rho (x^{(n)}) \in [0,1/2)$.) 
Finally, define 
\begin{equation} 
g(x)=t(Tx)-t(x)+h(\pi x) , 
\end{equation} 
thereby obtaining a positive continuous function on $X$ which takes only
finitely many values and also 
automatically satisfying (\ref{simple}).  Further, (\ref{simple}) implies
that $\phi[x,s]=[\pi x,0](s + t(x))$ is constant on fibers; therefore
it defines a continuous map $\flx \to
\fly$. 

Because $\phi$ satisfies (\ref{simple}), it commutes with the
$\Br$ actions. Because $\fly$ has a dense orbit (in fact all orbits
are dense),  $\phi$ is onto. Also, $\phi$ is one-to-one. For if
$\phi[x,s]=\phi[x',s']$, then $(x,s+t(x))\sim (x',s'+t(x'))$ (they
determine the same points in
$\fly$), so there is $d \in \Bz$ with $T^dx=x'$ and
$s'+t(x')=s+t(x)-h(x,d)$. But (\ref{simple}) implies that
\begin{equation}
h(x,d)-g(x,d)=t(x)-t(T^dx) ,
\end{equation}
so that 
\begin{equation*}
x'=T^dx, \quad s'=s-g(x,d) ,
\end{equation*}
and hence $(x,s)\sim (x',s')$ in $\flx$.

Since $t(x)=\rho (x)$ depends on the {\em entire sequence} $x$ of tile
types, $\phi [x,s] = [\pi x, 0](s+t(x))$ cannot be determined from a
finite window: the central coordinate of $\pi x$ can be, so that we
can determine from a finite window what tile in $\phi[x,s]$ is at the
origin, but we cannot tell exactly where in this tile to place the
origin without knowing the full sequence $x$ of tile types. 

\begin{rem}
The pair of maps $(T, \pi)$ in this example provides an action of $\Bz
^2$ on $X$, and (\ref{simple}), when expressed as
\begin{equation}
(t+g)(x) =t(Tx) + g(\pi x) ,
\end{equation}
is a sort of $\Bz ^2$ cocycle equation. See
\cite[Prop. 1.2]{KMS1}, \cite[p. 352]{KMS2}.
\end{rem}

\subsection{Example 2}\label{Ex2}
If one desires a simple map which is not a local code and, unlike the
preceding example, is not a topological conjugacy, it is easy to
modify Example \ref{Ex1} to accomplish that.
Let $(X,T)$ be the subshift of $\{0,1,2,3\}^{\Bz}$ determined by
coding $([0,1),R_{\alpha})$ according to entries to the four intervals
$[0,1/4)$, $[1/4,1/2)$, $[1/2,3/4)$, $[3/4,1)$. $(Y,S)$ is the coding of
$([0,1), R_{2\alpha})$ by the partition $\{[0,1/2),[1/2,1)\}$.
Now the factor map $\pi :(X,T) \to (Y,S)$ is determined by the 1-block
code that sends the symbols 0 and 2 to 0, and 1 and 3 to 1. (Under the
isomorphisms of $X$ and $Y$ with $[0,1)$, $\pi$ corresponds to $s \to
2s$ on $[0,1)$.) Then $h$ is defined as before except with $2\gamma$
replacing $\gamma$, and $t$ and $g$ are determined as before.
We arrive at a simple code which is not local, and which this time
is not one-to-one.

\subsection{Example 3}\label{Ex3}
With slightly more care in the choice of the parameters in Example
\ref{Ex1}, we can guarantee that {\em no local code exists} between
$\flx$ and $\fly$, even though there is a simple code $\phi: \flx \to \fly$.
In the definition of the ceiling function $h$ over $(Y,S)$, let the
division point $\gamma=1-\alpha$, and choose the heights $\eta_1$ and
$\eta_2$ so that $1,\alpha , \eta_1$, and $\eta_2$ are linearly
independent over $\Bz$. Then the ceiling function $g$ over $(X,S)$
takes values $\eta_1+\alpha$ (on $[0,1-\alpha )$) and
$\eta_2+\alpha-1$ (on $[1-\alpha ,1)$).

Suppose that $\phi : \flx \to \fly$ is a local code. Select $x \in X$
and let $[y,u]=\phi [x,0]$.
Now the central tile $y_0$ of $[y,u]$ and the position $u$ of the
origin in it are supposed to be completely determined by a finite
window in $[x,0]$, that is to say, by a finite block $C=x_{-r}\dots x_r$
of entries in $x$ (since the position of the origin in $x_0$ is fixed
at the left edge). Whenever $C$ appears in $x$, say the translated
tiling $[x,0]s$ agrees exactly with $[x,0]$ on its central range of
$2r+1$ tiles, we should see in the corresponding position of $[y,u]$
the same tile $y_0$ as before, with its left edge $u$ units to the
left of that of the central tile in $C$: 
\begin{equation}
(\pi_Y([y,u]s))_0=y_0, \quad \pi_{\Br}([y,u]s)=u .
\end{equation}
This forces a sum of lengths of tiles in $\flx$ to equal a sum of
lengths of tiles in $\fly$, which is impossible because of the linear
independence.

\subsection{Example 4}\label{Ex4}
Next we note that it is easy to construct local codes which are not
simple.
Let $(X,T)$ be a Sturmian symbolic dynamical system as in
Example \ref{Ex1}, with rotation $R_\alpha$ on $[0,1)$ coded by
entries to the interval $[0,1/2)$. Suppose that $\alpha > 1/2$, so
that the block $11$ does not appear in any $x \in X$. 
This time use the ceiling function
$g(x)\equiv 1$ on $X$, so that the associated tiling dynamical system
$\flx$ has
two tile types ($0$ and $1$), each of length $1$. Define a local code
onto another tiling system $\fly$ by splitting each tile labeled $1$
into two tiles, each of length $1/2$, both still labeled $1$. The
underlying symbolic dynamical system $(Y,S)$ is the primitive
(induced) transformation over $(X,S)$ formed by doubling $1$'s in all
sequences in $X$, and the ceiling function $h$ on $Y$ takes the value
$1$ on the cylinder over 0 (at time 0) and $1/2$ on the cylinder over
$1$. 

This local code is a topological conjugacy between the two tiling
dynamical systems, both of which are minimal. But $(Y,S)$ is
topologically weakly mixing \cite{K,PS} while $(X,T)$ has many
continuous eigenfunctions (every integer power of $\exp (2 \pi i\alpha)$ is an
eigenvalue), so there is no factor map from $(X,T)$ to $(Y,S)$ or vice
versa, and hence there is no simple code between $\flx$ and $\fly$.

\subsection{Example 5}\label{Ex5}
Finally, we construct a factor map between tiling dynamical systems
which is neither simple nor local, and in fact is such that no simple
or local code exists between the two tiling dynamical systems.

Let $(X,T)$ and $(Y,S)$ be as in Example \ref{Ex4}. $(X,T)$ is a
Sturmian symbolic dynamical system based on rotation modulo 1 by
$\alpha > 1/2$, coded via entries to an interval $[0,\beta)$ with
$\beta < 1- \alpha$ and $\beta \notin \Bz \alpha \mod 1$.  $(Y,S)$ is the
induced transformation over $(X,T)$ with a `second floor' over
$X_1=[1]=$ the cylinder set determined by fixing the central entry to
be $1$. If $X_0=[0]$ is the cylinder set consisting of those points
with central entry $0$, then $X=X_0 \cup X_1$, while
$Y=X_0 \cup X_1 \cup X_1'$, where $X_1'$ is a homeomorphic copy of
$X_1$ under a map $\theta: X_1 \to X_1'$.  Moreover,
\begin{equation}
Sy = \begin{cases}
     Ty \quad &\text{  if } y \in X_0\\
     \theta y &\text{  if  } y \in X_1\\
     T\theta^{-1}y &\text{  if  } y \in X_1' .
\end{cases}
\end{equation}
Define a continuous, not shift-commuting map $\pi : Y \to X$ which is
one-to-one on part of $Y$ and two-to-one on part by collapsing each
block $11$ to a single $1$; in terms of the induced transformation, 
\begin{equation}
\pi y = \begin{cases}
	y \quad &\text{  if  } y \in X_0 \cup X_1\\
	\theta ^{-1}y &\text{  if  } y \in X_1' ,
	\end{cases}
\end{equation}
so that
$\pi$ identifies $[1] \subset Y$ and the second floor over it and is
the identity on $[0]$.

This time the code will go in the direction $\phi : \fy \to
\fx$---we 
alter the choice of letters a little. As in Examples \ref{Ex1}
and  \ref{Ex3},
let $\gamma=1-\alpha$ and define $h$ on $X$ by
letting it take the value $\eta_1$ on $[0,\gamma )$ and $\eta_2$ on
$[\gamma ,1)$. 
We assume that $\eta_1$ and $\eta_2$ are positive and not too small,
say at least $5$.
Define $v_0: Y \to \Br$ by
\begin{equation}
v_0(y) = \rho (\pi y) =\begin{cases}
	\rho(y) \quad &\text{  if  } y \in X_0 \cup X_1\\
	\rho(\theta^{-1}y) &\text{  if  } y \in X_1' .
	\end{cases}
\end{equation}
Then $v_0(Sy)-v_0(y)$ is continuous on $Y$ and takes only finitely many
values ($0, \alpha$, and $-1+\alpha$). 
Define
\begin{equation*}
v(y)=v_0(Sy) .
\end{equation*}
Then $v(Sy) - v(y)$ still takes only the values $0, \alpha$, and $-1+
\alpha$, and it takes the value $\alpha$ on $X_1$.

Having defined $\pi, h$, and $v$, we must now define $g$ on $Y$ so as
to satisfy the relations (\ref{cohom}), which in the current setting take
the form
\begin{equation}\label{check}
\begin{split}
T^{n_X(\pi y,g(y)+v(y))}\pi y  &= T^{n_X(\pi Sy, v(Sy))}(\pi Sy) ,\\
v(Sy)-v(y)&=g(y) + h(\pi Sy,n_X(\pi Sy,v(Sy)))\\
          &\quad\qquad\qquad -h(\pi y, n_X(\pi y, g(y)+v(y))) .
\end{split}
\end{equation}
First note that if
\begin{equation}
m(y)=\chi_{X_0 \cup X_1'}(y) ,
\end{equation}
then
\begin{equation}
\pi(Sy)=T^{m(y)}(\pi y).
\end{equation}
Further, $n_X(x,v(y))=0$ for all $x,y$, since $0 \leq v(y) <1$ for all
$y$ and $h \geq 5$. Therefore
\begin{equation}
h(\pi Sy,n_X(\pi Sy,v(Sy)))=0 \quad\text{ for all } y,
\end{equation}
simplifying (\ref{check}).

Define
\begin{equation}
g(y)=v(Sy)-v(y)+h(\pi y,m(y)).
\end{equation}
Then $g$ is a continuous function on $Y$ that assumes finitely many
positive values (since $v(Sy)-v(y)= \alpha > 0$ on $X_1$, where
$m(y)=0$
so $h(\pi y,m(y))=0)$.
If we can verify that
\begin{equation}
n_X(\pi y, g(y)+v(y))=m(y) \quad\text{for all } y,
\end{equation}
then (\ref{check}) will follow.

Now $h(\pi y, m(y))$ is either $0$ or $h(\pi y)$, 
depending on whether
$m(y)=0$ or $1$. Thus
\begin{equation}
g(y)+v(y)=v(Sy)+h(\pi y,m(y)) ,
\end{equation}
and $v$ takes values in $[0,1)$, so the index $n_X(\pi
y,g(y)+v(y))$, which gives the value of $n$ for which
\begin{equation}
\sum_{k=0}^{n-1}h(T^k\pi y) \leq g(y)+v(y) < \sum_{k=0}^{n}h(T^k\pi y),
\end{equation}
must coincide with $m(y)$.

As in Example \ref{Ex4}, since there can be no factor map between the
underlying symbolic dynamical systems, for this example there can be no
simple code between the tiling dynamical systems. 
The tiles in $\fx$ have lengths $\eta_1$ and $\eta_2$, while those in
$\fy$ have lengths $\alpha,
-1+\alpha+\eta_1,\eta_2,\eta_2+\alpha ,$ and $-1+\alpha+\eta_2$.
If $1,\alpha,\eta_1$, and $\eta_2$ are chosen to be linearly
independent over $\Bz$, then
no nonzero sum of tile lengths
in $\fx$ can equal a sum of tile lengths in $\fy$, 
unless only tiles of length $\eta_2$ are used in both cases.
But since the systems $(X,T)$ and $(Y,S)$ are minimal, any finite
block $C=x_{-r}\dots x_r$ of entries in a sequence $x \in X$ will
appear separated by (arbitrarily long) blocks in which every symbol
appears. Thus the argument of Example \ref{Ex3} applies to show that
no local code can exist
between these tiling dynamical systems.

\bibliographystyle{amsplain}
\bibliography{tilecode}

\end{document}